\newif\ifforarxiv
\forarxivtrue

\ifforarxiv
\documentclass[11pt,a4]{amsart}
\usepackage{a4wide}
\else
\documentclass[reqno]{ijmart}
\fi

\usepackage[latin1]{inputenc}
\usepackage[T1]{fontenc}
\usepackage{lmodern}

\usepackage{graphicx,color,xy,latexsym,pstricks,verbatim}
\usepackage{amsmath,amssymb,amsfonts,amsbsy,amsthm}

\ifforarxiv
\usepackage[colorlinks=true,linkcolor=blue,citecolor=blue,bookmarks=false,pdftex]{hyperref}
\hypersetup{
  pdfauthor   = {J.-M. Couveignes and R. Lercier},
  pdftitle    = {Fast construction of irreducible polynomials over finite fields},
  pdfsubject  = {Algorithmic Number Theory},
  pdfkeywords = {11G20, 11Y16, 14H45, 68Q25}
}
\fi

\definecolor{LightGrey}{rgb}{.85,.85,.85}
\definecolor{DarkGrey}{rgb}{.5,.5,.5}

\definecolor{Blue}{rgb}{.0,.0,0.9}
\definecolor{LightBlue1}{rgb}{.2,.4,0.9}
\definecolor{LightBlue2}{rgb}{.3,.5,0.9}
\definecolor{LightBlue3}{rgb}{.4,.6,0.9}
\definecolor{LightBlue4}{rgb}{.5,.7,.9}
\definecolor{LightBlue5}{rgb}{.6,.8,.9}
\definecolor{LightBlue6}{rgb}{.7,.9,.9}

\definecolor{Red}{rgb}{.9,.0,.0}
\definecolor{LightRed1}{rgb}{0.9,.2,.4}
\definecolor{LightRed2}{rgb}{0.9,.3,.5}
\definecolor{LightRed3}{rgb}{0.9,.4,.6}
\definecolor{LightRed4}{rgb}{.9,.5,.7}
\definecolor{LightRed5}{rgb}{.9,.6,.8}
\definecolor{LightRed6}{rgb}{.9,.7,.9}

\definecolor{Grey}{rgb}{.5,.5,.5}

\definecolor{Blue}{rgb}{.0,.0,0.9}
\definecolor{LightBlue1}{rgb}{.2,.4,0.9}
\definecolor{LightBlue2}{rgb}{.3,.5,0.9}
\definecolor{LightBlue3}{rgb}{.4,.6,0.9}
\definecolor{LightBlue4}{rgb}{.5,.7,.9}
\definecolor{LightBlue5}{rgb}{.6,.8,.9}
\definecolor{LightBlue6}{rgb}{.7,.9,.9}

\definecolor{Red}{rgb}{.9,.0,.0}
\definecolor{LightRed1}{rgb}{0.9,.2,.4}
\definecolor{LightRed2}{rgb}{0.9,.3,.5}
\definecolor{LightRed3}{rgb}{0.9,.4,.6}
\definecolor{LightRed4}{rgb}{.9,.5,.7}
\definecolor{LightRed5}{rgb}{.9,.6,.8}
\definecolor{LightRed6}{rgb}{.9,.7,.9}

\xyoption{all}
\usepackage[english]{babel}

\newcounter{noalgo}[section]
\newdimen\indentalgo
\newdimen\indentalgodec\indentalgo=0.0mm\indentalgodec=10mm

\newcommand{\If}{\advance\indentalgo by \indentalgodec {\bf if }}

\newcommand{\For}{\global\advance\indentalgo by \indentalgodec {\bf for }}

\newcommand{\Endindent}{\global\advance\indentalgo by -\indentalgodec}

\newdimen\decalage \decalage=0.5cm
\newcounter{algo} 

\let\set\mathbb
\def\<<{\leavevmode
  \raise0.28ex\hbox{$\scriptscriptstyle\langle\!\langle$}\nobreak
  \hskip -.6pt plus.3pt minus.2pt\,}
\def\>>{\,\nobreak\hskip -.6pt plus.3pt minus.2pt
  \raise0.28ex\hbox{$\scriptscriptstyle\rangle\!\rangle$}}

\def\Aut{\mathop{\rm{Aut}}\nolimits }

\def\Ker{\mathop{\rm{Ker}}\nolimits }
\def\Gal{\mathop{\rm{Gal}}\nolimits }

\def\End{\mathop{\rm End }}
\def\FF{{\set F}}

\def\FQ{{\FF _Q}}
\def\Fp{{\FF _p}}

\def\Fq{{\FF _q}}

\def\Fqs{{{\FF ^*_q}}}
\def\FQs{{{\FF ^*_Q}}}

\def\bK{{\bf K  }}
\def\bM{{\bf M  }}
\def\bA{{\bf A  }}
\def\bS{{\bf S  }}
\def\bL{{\bf L  }}

\def\bG{{\bf  G }}
\def\Gm{{{\bf  G }_m}}

\def\ZZ{{\set Z}}

\def\cA{{\mathcal A}}

\def\cE{{\mathcal E}}

\def\cG{{\mathcal G}}

\def\cI{{\mathcal I}}

\def\cM{{\mathcal M}}

\def\cO{{\mathcal O}}

\def\cW{{\mathcal W}}

\def\dgot{{\mathfrak  d}}
\def\lgot{{\mathfrak l}}

\newtheorem{lemma}{Lemma}
\newtheorem{corollary}{Corollary}

\newtheorem{theorem}{Theorem}

\usepackage{enumitem}

\begin{document}

\ifforarxiv

\author{Jean-Marc Couveignes}
\address{INRIA Bordeaux Sud-Ouest \& Universit\'e de Toulouse.}
\address{Universit\'e de Toulouse le Mirail,
  D{\'e}partement de Math{\'e}matiques et Informatique,\linebreak
  5 all{\'e}es Antonio Machado,
  F-31058 Toulouse c{\'e}dex 9,
  France.
}
\email{jean-marc.couveignes@math.univ-toulouse.fr}

\author{Reynald Lercier} 
\address{
  DGA,
  La Roche Marguerite,
  F-35174 Bruz Cedex,
  France.}
\address{
  Institut de Recherche Math\'ematique de Rennes,
  Universit\'e de Rennes 1,
  Campus de Beaulieu,
  F-35042 Rennes,
  France.}
\email{reynald.lercier@m4x.org}

\date{July 28$^{\textrm{th}}$, 2011}
\else

\author{Jean-Marc Couveignes}
\address{INRIA Bordeaux Sud-Ouest \& Universit\'e de Toulouse.\linebreak}
\address{Universit\'e de  Toulouse le Mirail,
  D{\'e}partement de Math{\'e}matiques et Informatique,\linebreak
  5 all{\'e}es Antonio Machado,
  F-31058 Toulouse c{\'e}dex 9,
  France.
}
\email{jean-marc.couveignes@math.univ-toulouse.fr}

\author{Reynald Lercier} 
\address{
  DGA,
  La Roche Marguerite,
  F-35174 Bruz Cedex,
  France.\linebreak
}
\address{
  Institut de Recherche Math\'ematique de Rennes,
  Universit\'e de Rennes 1,\linebreak
  Campus de Beaulieu,
  F-35042 Rennes,
  France.}
\email{reynald.lercier@m4x.org}

\fi

\title[Irreducible polynomials over finite fields]{Fast construction of irreducible polynomials over finite fields}

\bibliographystyle{plain}

\begin{abstract}
  We present a randomized algorithm that on input a finite field $\bK$ with
  $q$ elements and a positive integer $d$ outputs a degree $d$ irreducible
  polynomial in $\bK[x]$.  The running time is $d^{1+\varepsilon(d)} \times
  (\log q)^{5+\varepsilon(q)}$ elementary operations. The function
  $\varepsilon$ in this expression is a real positive function belonging to
  the class $o(1)$, especially, the complexity is quasi-linear in the
  degree $d$. Once given such 
an  irreducible polynomial of degree $d$, we can 
   compute  random irreducible polynomials of degree $d$ at the expense
  of  $d^{1+\varepsilon(d)} \times (\log q)^{1+\varepsilon(q)}$ elementary
  operations only.
\end{abstract}

\maketitle

%\tableofcontents

\section{Introduction}

This article deals with the following problem: given a prime integer $p$, a
power $q=p^w$ of $p$, a finite field $\bK$ with $q$ elements, and a positive
integer $d$, find a degree $d$ irreducible polynomial in $\bK[x]$.  We assume
that the finite field $\bK$ is given as a quotient $(\ZZ/p\ZZ)[z]/(m(z))$
where $m(z)$ is a degree $w$ monic irreducible polynomial in
$(\ZZ/p\ZZ)[z]$. We assume furthermore that polynomials are given in a dense
representation. The complexity of algorithms will be evaluated in terms of the
number of necessary elementary operations. Additions, subtractions and
comparisons in $\bK$ require a constant times $\log q$ elementary
operations. Multiplication and division require $(\log q)\times (\log\log
q)^{1+\varepsilon(q)}$ elementary operations. Note that in this paper, the
notation $\varepsilon(x)$ stands for a real positive function of the parameter
$x$ alone, belonging to the class $o(1)$.

It has been proven by Adleman and Lenstra~\cite{AL} that under the generalized
Riemann hypothesis there exists an algorithm that constructs a degree $d$
irreducible polynomial over $\bK$ in {\it deterministic} polynomial time in
$d$ and $\log q$. There is no known unconditional proof of this result.  The
main algorithms in this paper are \textit{Las Vegas probabilistic}.  The
behavior of a Las Vegas algorithm depends on the input of course, but also on
the result of some random choices. One has to flip coins. 
A Las Vegas  algorithm either stops with the correct result 
or informs that it
failed.
 The running time of the algorithm is bounded from above in terms of the size
 of the input only (this upper bound should not depend on the random
 choices).  For {\it each} input, one asks that the probability that the
 algorithm succeeds  is $\geqslant 1/2$.
See Papadimitriou's book~\cite{papadimitriou} for a formal
definition of the main complexity classes.

A classical probabilistic approach to finding irreducible polynomials consists
in first choosing a random polynomial of degree $d$ and then testing for its
irreducibility.
The probability that a polynomial of degree $d$ be irreducible is greater than
or equal to $1/(2d)$. See Lidl and Niederreiter~\cite[Ex. 3.26 and 3.27, page
142]{Lidl} and Lemma~\ref{lemma:lidl2} of Section~\ref{section:compter} below.
In order to check whether a polynomial $f(x)$ is irreducible, we may use
Ben-Or's irreducibility test~\cite{BO}.  This test has maximal complexity
$(\log q )^{2+\varepsilon(q)}\times d^{2+\varepsilon(d)}$ elementary
operations while its average complexity is $(\log q)^{2+\varepsilon(q)} \times
d^{1+\varepsilon(d)}$ elementary operations according to Panario and
Richmond~\cite{PR}.  The average complexity of finding an irreducible
polynomial with this method is thus $d^{2+\varepsilon (d)}\times (\log
q)^{2+\varepsilon(q)}$ elementary operations.  All the known algorithms have a
quadratic factor at least in $d$. A survey can be found in the work of
Lenstra~\cite{L1,L2} and Shoup~\cite[section 1.2]{Shoup93}.  It seems
difficult to improve on these existing methods as long as we use an
irreducibility test.

So we are driven to consider very particular polynomials.
For example, Adleman and Lenstra~\cite{AL} construct irreducible polynomials
in this way.  Their method uses Gauss periods that are relative
traces of roots of unity. In Section~\ref{section:qm1},
we recall how efficient known methods can be for very special values of the
degree $d$. We reach quasi-linear complexity in $d$ when $d=\ell^\delta$ is a
power of a prime divisor $\ell$ of $p(q-1)$.
Section~\ref{section:general} explains how to construct a degree $d_1d_2$
irreducible polynomial once given two irreducible polynomials of coprime
degrees $d_1$ and $d_2$.  We explain in Sections~\ref{section:fiber}
and~\ref{section:basis} how to construct irreducible polynomials using
isogenies between elliptic curves.  Thanks to this new construction, we reach
quasi-linear complexity in $d$ when $d=\ell^\delta$ is a power
of a prime   $\ell$  and $\ell$ does not divide  $p(q-1)$.
Putting everything together, we obtain a probabilistic algorithm that finds
degree $d$ irreducible polynomials in $\bK[x]$ in quasi-linear time in $d$,
without any restriction on $d$ nor $q$.  Our constructions are summarized in
Section~\ref{section:compter}.  In Section~\ref{sec:prel-results-about}, we
state several useful preliminary results about finite fields, polynomials and
elliptic curves.
\medskip

\noindent{\it Remark.} One may wonder if the algorithms and complexity
estimates in this paper are still valid when the base field is not presented
as a quotient of the form $(\ZZ/p\ZZ)[z]/(m(z))$.  Following~\cite[Section 1]{L2}, one may
assume for example that elements in $\bK$ are represented as vectors in
$(\ZZ/p\ZZ)^w$. Assume we are given the vector corresponding to the unit
element $1$. Assume also we are given a black box or an algorithm that
computes multiplications and divisions of elements in $\bK$.  In this
situation, before applying the algorithms presented in this paper, we should
first construct an isomorphism between the given $\bK$ and a quotient ring of
the form $(\ZZ/p\ZZ)[z]/(m(z))$.  To this end, we first look for a generator
$\tau$ of the $(\ZZ/p\ZZ)$-algebra $\bK$. We pick a random element $\tau$ in
$\bK$. The probability that $\tau$ generates $\bK$ over $\ZZ/p\ZZ$ is at least
$1/2$ according to Lemma~\ref{lemma:lidl2} of Section~\ref{section:compter}.
We compute the powers $\tau^k$ for $0\leqslant k \leqslant w$. These are $w+1$
vectors of length $w$. We compute the kernel of the corresponding matrix in
$\cM_{w\times (w+1)}(\ZZ/p\ZZ)$.  If the dimension of this kernel is bigger
than $1$, then $\tau$ is not a generator, so we pick a different $\tau$ and
start again. If the kernel has dimension $1$, we obtain the minimal polynomial
$m(z)\in (\ZZ/p\ZZ)[z]$ of $\tau$, and an explicit isomorphism $\kappa$ from
$\tilde \bK=(\ZZ/p\ZZ)[z]/(m(z))$ onto $\bK$. All this requires $O(w)$
operations in $\bK$ and $O(w^3)$ operations in $\ZZ/p\ZZ$. Given any degree
$d$ irreducible polynomial $\tilde f(x)$ in $\tilde \bK[x]$, we deduce an
irreducible polynomial in $\bK[x]$ by applying the isomorphism $\kappa$ to
every coefficient in $\tilde f(x)$. This requires $dw^2$ operations in
$\ZZ/p\ZZ$.  So our algorithms and complexity estimates remain valid in that
case, as long as elementary operations in $\bK$ can be computed in time $(\log
q)^{4+\varepsilon(q)}$ elementary operations. This includes all the reasonable
known models for finite fields, including normal bases, explicit
data~\cite{L2} and towers of extensions.
\medskip

\noindent{\it Notations.}  If $\bK$ is a field with characteristic $p$ and $q$
is a power of $p$, we denote by $\Phi_q: \bK \rightarrow \bK$ the morphism
which raises to the $q$-th power.  If $\bG$ is an algebraic group over 
a field with $q$ elements,
we denote by $\varphi_{\bG} : \bG\rightarrow \bG^{(q)}$, the Frobenius
morphism.\smallskip 

The book~\cite{liu} by Liu provides a nice introduction to abstract
algebraic geometry.  A good monograph on elliptic curve is~\cite{silv}.
\medskip

\noindent {\it Acknowledgments.} We thank K. Kedlaya for pointing out his
joint work with Umans~\cite{UK} to us, and H. Lenstra for explaining to us how
to save a $\log q$ factor in the complexity using~\cite{Howe}.
We thank C{\'e}cile Dartyge, Guillaume Hanrot,  Gerald Tenenbaum and Jie Wu
for pointing out Iwaniec's result on Jacobsthal's problem  to us.

This work was partially supported by the French ``Agence Nationale de la
Recherche'' under the \textsc{anr-07-blan-0248 algol} and by the French
``D{\'e}l{\'e}gation G{\'e}n{\'e}rale pour l'Armement''.

\section{Basic constructions}\label{section:qm1}

In this section, $\bK$ is a finite field with $q=p^w$ elements and $\Omega$ is
an algebraic closure of $\bK$. For every positive integer $k$, we denote by
$\FF_{p^k}$ the unique subfield of $\Omega$ with $p^k$ elements. We explain
how to quickly construct a degree $d$ irreducible polynomial when $d$ is a
prime power $\ell^\delta$ and $\ell$ divides $p(q-1)$.  All the constructions
in this section are known, but deserve to be quickly
surveyed. Section~\ref{subsection:AS} deals with the case
$\ell=p$. Section~\ref{subsection:rad} deals with the case when $\ell$ is a
prime divisor of $(q-1)$. Section \ref{subsection:special} is concerned with
the special case $\ell=2$ and $q$ odd. In Section~\ref{subsection:gauss}, we
illustrate on a simple example how to use Kummer theory together with descent
when the roots of unity lie in a non-trivial extension of the base field.
Although the results in Section \ref{subsection:gauss} are not necessary to
prove our main theorems in Section~\ref{section:compter}, several ideas at
work in this section play a decisive role later in the slightly more advanced
context of Section \ref{section:basis}.

\subsection{Artin-Schreier towers}\label{subsection:AS}

In this section, we are given a $p$-th power $d=p^\delta$ and we want to
construct a degree $d$ irreducible polynomial in $\bK[x]$.  We use a
construction of Lenstra and de Smit~\cite{LdS} in that case.  
If $k$ and $l$ are two positive integers such that $l$ divides $k$, we define
the polynomial
$T_{l,k}(x)=x+x^{p^l}+x^{p^{2l}}+\cdots+x^{p^{(\frac{k}{l}-1)l}}$.  For every
positive integer $k$, we denote by $\cA_k\subset \Omega$ the subset consisting
of all scalars $a\in \Omega$ such that the three following conditions hold
true:
\begin{itemize}[itemsep=0pt,parsep=0pt,partopsep=0pt,topsep=0pt]
\item $a$ generates $\FF_{p^k}$ over $\Fp$, \emph{i.e.} $\Fp(a)=\FF_{p^k}$\,,
\item $a$ has non-zero absolute trace, \textit{i.e.} $T_{1,k}(a)\not =0$\,,
\item $a^{-1}$ has non-zero absolute trace, \textit{i.e.} $T_{1,k}(a^{-1})\not
=0$.
\end{itemize}
We set $I(X)=({X^p-1})/({\sum_{1\leqslant i\leqslant p-1}X^i})$.  This
rational fraction induces a  $p$ to $1$ surjective map
\begin{displaymath} 
  I : \Omega-\Fp \rightarrow \Omega-\{0\}.
\end{displaymath}
We check that $I^{-1}(\cA_k)\subset \cA_{pk}$ for every $k\geqslant 1$.
Indeed, if $a\in \cA_k$ and if $I(b)=a$ then $b\not =1$ and
\begin{displaymath} 
  \frac{1}{(1-b)^p}-\frac{1}{1-b}=
  \frac{b^p-b}{(b-1)^{p+1}}= \frac{b+\cdots+b^{p-1}}{b^p-1}=a^{-1}.
\end{displaymath}
So $1/(1-b)$ is a root of the separable polynomial $x^p-x=a^{-1}$.  
This polynomial is irreducible over $\FF_{p^k}[x]$ because the absolute trace
of $a^{-1}$ is non-zero. So $\Fp(b)=\FF_{p^{pk}}$. Further, $b$ is a root of
the polynomial $x^{p}-a(x^{p-1}+\cdots+x)-1$. So the trace $T_{k,pk}(b)$ of
$b$ relative to the extension $\FF_{p^{pk}}/ \FF_{p^{k}}$ is $a$. As a
consequence the absolute trace of $b$ is
$T_{1,pk}(b)=T_{1,k}(T_{k,pk}(b))=T_{1,k}(a)$, the absolute trace of $a$, and
it is non-zero.  Now $b^{-1}$ is a root of the reversed polynomial
$x^{p}+a(x^{p-1}+\cdots+x)-1$. So the trace of $b^{-1}$ relative to the
extension $\FF_{p^{pk}}/ \FF_{p^{k}}$ is $-a$. As a consequence, the absolute
trace of $b^{-1}$ is the opposite of the absolute trace of $a$, and it is
non-zero.

Since $\cA_1=\Fp-\{0\}$, we deduce that $\#\cA_{p^k}\geqslant (p-1)p^{k}$.  In
particular the fiber above $1$ of the iterated rational fraction
$I^{(\delta)}$ is irreducible over $\Fp$. 
See \cite[Section 3.1]{liu} for the definition of the fiber of a morphism over
a point, and \cite[Section 2.4]{liu} for the definition of an irreducible scheme.
If $w$ is prime to $p$, then this
fiber remains irreducible over $\bK=\Fq$. In general, we factor the degree $w$
of $\Fq/\Fp$ as $w=p^ew'$ where $w'$ is prime to $p$.  We first look for an
element $a\in \cA_{p^e}\subset \Fq$.  Using the remarks above we can find such
an $a$ by solving $e$ Artin-Schreier equations with coefficients in $\Fq$. To
this end, we write down the matrix of the $\Fp$-linear map $x\mapsto x^p-x$ in
the $\Fp$-basis $(1,z,\ldots,z^{w-1})$ of $\bK=(\ZZ/p\ZZ)[z]/(m(z))$.  We then
solve the $e$ corresponding $\Fp$-linear systems of dimension $w$.
Altogether, finding $a$ requires a constant times $w \times \log p$ operations
in $\bK$ and a constant times $ew^3$ operations in $\Fp$. Since $w=O(\log q)$
and $e=O(\log w)=O(\log \log q)$ we end up with a complexity of $(\log
q)^{4+\varepsilon(q)}$ elementary operations.

The fiber $I^{-\delta}(a)$ is a degree $p^\delta$ irreducible divisor over
$\FF_{p^{p^e}}$. It remains irreducible over $\bK=\Fq$. It remains to
compute the annihilating polynomial of this fiber.  We compute the iterated
rational fraction $I^\delta(x)={N(x)}/{D(x)}$. Composition of polynomials and
power series can be computed in quasi-linear time
i.e. $d^{1+\varepsilon(d)}\times (\log q)^{1+\varepsilon(q)}$ elementary
operations, using recent results by Umans and Kedlaya~\cite{Umans08,UK} (see
Corollary~\ref{cor:compof} in Section~\ref{subsection:comp}).  An older
algorithm due to Brent and Kung has exponent $({\omega+1})/{2}+\varepsilon(d)$
where $\omega$ is the exponent in matrix multiplication. So we can compute
$N(x)$ and $D(x)$ at the expense of $p^{\delta(1+\varepsilon(p^\delta))}\times
(\log q)^{1+\varepsilon(q)}$, that is $d^{1+\varepsilon(d)}\times (\log
q)^{1+\varepsilon(q)}$ elementary operations. The polynomial $f(x)=N(x)-aD(x)$
is an irreducible degree $d$ polynomial in $\bK[x]$.

We thus have proven the following lemma.

\begin{lemma}[Artin-Schreier extensions]\label{lemma:ASE}
  There exists a deterministic algorithm that on input a finite field
  $\bK=(\ZZ/p\ZZ)[z]/(m(z))$ with cardinality $q=p^w$ and a positive integer
  $\delta$ computes an irreducible degree $d=p^\delta$ polynomial in $\bK[x]$
  at the expense of $(\log q)^{4+\varepsilon(q)}+ d^{1+\varepsilon(d)}\times
  (\log q)^{1+\varepsilon(q)}$ elementary operations.
\end{lemma}
\medskip

\noindent {\it Example.} We take $p=2$, $q=4$, $\delta=2$ and $d=4$.  We
assume $\bK=\FF_2[z]/(z^2+z+1)$, so $e=1$. We know that $1\in \cA_1$. We set
$a=z\bmod z^2+z+1$ and check that $I(a)=1$, so $a\in \cA_2$. We compute
$I(I(x))=(x^4+x^2+1)/(x^3+x)$ and set $f(x)= x^4+x^2+1- a(x^3+x)$. This is an
irreducible polynomial in $\bK[x]$.

\subsection{Radicial extensions}\label{subsection:rad}

In this section, $\ell$ is a prime dividing $q-1$.  Let $d=\ell^\delta$ for
some positive integer $\delta$.  In the special case 
 $\ell=2$ we further ask
that $\ell^2=4$ divide $q-1$.  We want to construct a degree $d$ irreducible
polynomial in $\bK[x]$.  This is a very classical case. 

We write
$q-1=\ell^e\ell'$ where $\ell'$ is prime to $\ell$.  We first look for a
generator $a$ of the $\ell$-Sylow subgroup of $\Fqs$. To find such a
generator, we pick a random element in $\Fqs$ and raise it to the power
$\ell'$. Call $a$ the result. Check that $a^{\ell^{e-1}}\not =1$.  If this is
not the case, start again. The probability of success is $1-1/\ell$. The
average complexity of finding such an $a$ is $O(\log q)$ operations in
$\Fq$. The polynomial $f(x)=x^{d}-a$ is irreducible in $\Fq[x]$. This is well
known but we try to prove it in a way that will be easily adapted to a more
general context later.

The $\ell^{\delta+e}$-torsion $\Gm[\ell^{\delta+e}]$ of the multiplicative
group $\Gm / \Fq$ is isomorphic to $(\ZZ/\ell^{\delta+e}\ZZ,+)$ and the Frobenius
endomorphism $\varphi_{\Gm} : \Gm\rightarrow \Gm$ acts on it as multiplication
by $q$. The order of $q=1+\ell'\ell^{e}$ in $(\ZZ/\ell^{e+\delta}\ZZ)^*$ is
$\ell^\delta=d$. So the Frobenius $\Phi_q$ acts transitively on the roots of
$f(x)$.
\medskip

\noindent {\it Example.} We take $p=5$, $q=5$, $\ell=2$, $\delta=3$ and
$d=8$. We check that $4$ divides $p-1$. In particular $e=2$ and $\ell'=1$.
The class $a=2\bmod 5$ generates the $2$-Sylow subgroup of
$(\ZZ/5\ZZ)^*$. Indeed $2^4=1\bmod 5$ and $2^2=-1\bmod 5$.  We set
$f(x)=x^8-2$.

\subsection{A special case}\label{subsection:special}

In this section, we assume that $p$ is odd, $\ell=2$ and $d=2^\delta$ for some
positive $\delta$. We need to adapt the methods of
Section~\ref{subsection:rad} in that special case because the group of units
in $\ZZ/d\ZZ$ that are congruent to $1$ modulo $\ell$ is no longer cyclic when
$\ell=2$ and $\delta>2$.  We want to construct a degree $d$ irreducible
polynomial in $\bK[x]$.  This time we assume that $2^2$ does not divide
$q-1$. So $q$ is congruent to $3$ modulo $4$.  We set $Q=q^2$ and observe that
$4$ divides $Q-1$.

We first look for a generator $c$ of $\FQ$ over $\bK=\Fq$.  For example we
take $c$ a root of the polynomial $y^2-r$ where $r$ is not a square in $\bK$
(we can take $r=-1$). If $\delta=1$ we are done. Assume now $\delta\geqslant
2$.  We write $Q-1=2^e\ell'$ where $\ell'$ is prime to $2$.  We find a
generator $a$ of the $2$-Sylow subgroup of $\FQs$.  The polynomial
$F(x)=x^{d/2}-a$ is irreducible in $\FQ[x]$. It remains to derive from
$F(x)$ an irreducible polynomial $f(x)$ of degree $d$ in $\bK[x]$. We call
$\bar a =\Phi_q(a)=a^q$, the conjugate of $a$ over $\Fq$. We can compute it at
the expense of $O(\log q)$ operations in $\bK$. It is clear that $\bar a \not
= a$ because the order of $a$ is divisible by $4$ and there is no point of
order $4$ in $\Gm(\Fq)$. The polynomial $f(x)=(x^{d/2}-a)(x^{d/2}-\bar a)$ has
coefficients in $\bK$.  It is irreducible over $\bK$. Indeed, any root $b$ of
$x^{d/2}-a$ is also a root of $f(x)$.  The field $\Fq(b)$ generated by $b$
over $\Fq$ contains $a$ and it has degree $d/2$ over $\Fq(a)=\FQ$ because
$F(x)$ is irreducible in $\FQ[x]$. So $f(x)$ is irreducible in $\bK[x]$.
\medskip

Sections~\ref{subsection:rad} and~\ref{subsection:special} prove the following
lemma.
\begin{lemma}[Kummer extensions]\label{lemma:KE}
  There exists a probabilistic (Las Vegas) algorithm that on input a finite
  field $\bK=(\ZZ/p\ZZ)[z]/(m(z))$ with cardinality $q=p^w$, a prime integer
  $\ell$ dividing $q-1$, and a positive integer $\delta$, computes an
  irreducible degree $d=\ell^\delta$ polynomial in $\bK[x]$ at the expense of
  $(\log q)^{2+\varepsilon(q)}+ d \log q$ elementary operations.
\end{lemma}
\medskip

\noindent {\it Example.} We take $p=7$, $q=7$, $\ell=2$, $\delta=3$ and
$d=8$. Since $4$ does not divide $q-1$, we set $Q=q^2=49$. We factor
$49-1=2^4\times 3$ so $e=4$ and $\ell'=3$. We check that $r=3\bmod 7$ is not a
square in $\FF_7$. So we set $c=y\bmod y^2-3\in \FF_7[y]/(y^2-3)$. We set
$a=(1+c)^3=3-c$ and check $a^{16}=1$ and $a^8=-1$. We set $F(x)=x^4-a$. We
compute $\bar a =a^7=3+c$.  We set $f(x)=(x^4-a)(x^4-\bar a)=x^8+x^4-1$. This
is an irreducible polynomial in $\FF_7[x]$.

\subsection{Descent}\label{subsection:gauss}

In this section, we recall how to use Kummer theory when roots of unity are
missing. We do not hope to find a quasi-optimal algorithm that way. But
several important algorithmic questions arise naturally in this context.

We assume $\ell=3$ and $d=3^\delta$ and $p=q\not =3$.  We assume that $3$ does
not divide $q-1$. So $q$ is congruent to $2$ modulo $3$, and we cannot apply
the method in Section~\ref{subsection:rad}.  We experiment in this simple
context an idea that will be decisive in Section~\ref{section:basis}. We base
change to a small auxiliary extension.  We set $Q=q^2$ and observe that $3$
divides $Q-1$. We shall deal with the field $\FQ$ with $Q$ elements. We note
that this idea is valid for any prime $\ell$, but the degree of the auxiliary
extension $\FQ/\Fq$ might be quite large (up to $\ell-1$) for a general
$\ell$.

We first need to build a computational model for the field $\FQ$. For example
we pick a degree $2$ irreducible polynomial $y^2-r_1y+r_2$ in $\bK[x]$ and set
$\bL=\bK[y]/(y^2-r_1y+r_2)$.  We set $c=y\bmod y^2-r_1y+r_2$.  We write
$Q-1=3^e\ell'$ where $\ell'$ is prime to $3$.  We find a generator $a$ of the
$3$-Sylow subgroup of $\bL^*$.  The polynomial $F(x)=x^{d}-a$ is irreducible
in $\bL[x]$. It remains to derive from $F(x)$ an irreducible polynomial
$f(x)$ of degree $d$ in $\bK[x]$.

Let $b = x \bmod F(x)$. This is a root of $F(x)$ in $\bL[x]/(F(x))$. The
latter field has $q^{2d}$ elements.  Recall $\Phi_q$ is the map which
raises to the $q$-th power. We have $\Phi_Q= \Phi_q^2$.  For any $\alpha$ in
$\bL[x]/(F(x))$, we set $\Sigma_1(\alpha)= \alpha+\Phi_q^d(\alpha)$ and
$\Sigma_2(\alpha)=\alpha\times \Phi_q^d(\alpha)$.

  \begin{displaymath}
    \xymatrixrowsep{0.1in}
    \xymatrixcolsep{0.1in}
    \xymatrix{ \bL[x]/(F(x))\simeq \FF_{q^{2d}} \ar@{-}[rd]
      \ar@{-}[dd] &\\ &\bK(\Sigma_k(b))\simeq \FF_{q^{d}} \ar@{-}[dd]\\ \bL =
      \bK[y]/(y^2-r_1y+r_2)\simeq \FF_{q^2} \ar@{-}[rd] &\\ & \bK\simeq \FF_q \\ }
  \end{displaymath}

Since $d$ is a {\it prime power}, at least one among $\Sigma_1(b)$ and
$\Sigma_2(b)$ generates an extension of degree $d$ of $\Fq$ (see
Lemma~\ref{lemma:subfield} of Section~\ref{section:gene}).
In other words, there exists a $k\in \{1,2\}$ such that the polynomial
\begin{displaymath}
  f(x) = \prod_{0 \leqslant l < d}(x - \Phi_q^l(\Sigma_k(b)))
\end{displaymath} 
is irreducible of degree $d$ in $\bK[x]$\,.
Three questions now worry us.
\begin{enumerate}[itemsep=0pt,parsep=0pt,partopsep=0pt,topsep=0pt]
\item How to compute $\Sigma_k(b)$ for $k\in \{1,2\}$~?
\item How to find the good integer $k$~?
\item How to compute $f(x)$ starting from $F(x)$~?
\end{enumerate}
\medskip

Question $1$ boils down to asking how to compute $\Phi_q^d(b)$. A first method
would be to compute $\Phi_q^d(b)$ as $b^{q^{d}}$ at the expense of a constant
times $d\log q$ operations in $\bL[x]/(F(x))$. This would require a constant
times $\log q \times d^{2+\varepsilon(d)}$ operations in $\bK$. This is too
much for us.

Instead of that, we should remind ourselves of the geometric origin of the
polynomial $F(x)$. Indeed, $b$ lies in $\Gm[3^{e+\delta}]$.  We write
$q^d=R\bmod 3^{e+\delta}$ where $0\leqslant R<3^{e+\delta}\leqslant Qd$. Then
$\Phi_q^d(b)=b^R$ can be computed at the expense of a constant times $\log
R\le  \log Q+\log d$ operations in $\bL[x]/(F(x))$. This requires a constant times
$\log q \times d^{1+\varepsilon(d)}$ operations in $\bK$.
\medskip

Question $2$ can be solved by comparing $\Sigma_1(b)$ and its conjugate by
$\Phi_q^{3^{\delta-1}}$. We have
\begin{displaymath}
  \Phi_q^{3^{\delta-1}}(\Sigma_1(b))=\Sigma_1(\Phi_q^{3^{\delta-1}}(b))=\Phi_q^{3^{\delta-1}}(b)
  +\Phi_q^{3^\delta+ 3^{\delta-1}}(b).
\end{displaymath} 
Each of the two terms in the above sum can be computed as explained in the
answer of Question~$1$. Since $\Sigma_2(b)=1$ here, we already know that
$\Sigma_1(b)$ is the good candidate. But we keep the more naive approach in
mind.  \medskip

Question $3$ is related to the following problem: we are given $\Sigma_k(b)$
for some $k\in\{1,2\}$. We know that $\Sigma_k(b)$ generates the degree $d$
extension of $\bK$ inside $\bL[x]/(F(x))$.  Therefore its minimal polynomial
$f(x)$ in the latter extension has coefficients in $\bK$.  We want to compute
this degree $d$ polynomial in $\bK[x]$.  One can apply a general algorithm for
this task, such as the one given by Kedlaya and Umans (\cite{Umans08,UK} and
Theorem~\ref{th:mini} below).  They show that it is possible to compute this
minimal polynomial at the expense of $d^ {1+\varepsilon(d)}\times (\log
Q)^{1+\varepsilon(q)}$ elementary operations. Thus, the complexity is
quasi-linear in the degree $d$.  \medskip

\noindent {\it Example.} We take $p=q=5$, $\ell =3$, $\delta=2$, $d=9$. So
$Q=25$, $Q-1=3\times 8$, $e=1$ and $\ell'=8$.  We check that $r=2\bmod 5$ is
not a square. We set $c=y\bmod y^2-2 \in \FF_5[y]/(y^2-2)$. We compute
$a=(1+c)^8=2+3c$.  We check $a^3=1$ and $a\not =1$. We set $F(x)=x^9-a$ and
$b=x\bmod F(x)$. We need to compute the conjugate of $b$ above
$\FF_{5^{9}}$. This is $b^{5^9}$. Recall $b$ lies in $\Gm[27]$. So we don't
raise $b$ to the power $5^9$ brutally.  We rather compute $5^9=1953125=-1\bmod
27$.  So $\Phi_{5^9}(b)=1/b=2(y+1)x^8\bmod(x^9-2-3y,y^2-2,5)$. The product
$\Sigma_2(b)=1$ is not the good candidate. So we compute the minimal
polynomial of $\Sigma_1(b)=b+1/b$ and find $f(x)=x^9 + x^7 + 2x^5 + 4x + 1\in
\FF_5[x]$.

\section{Preliminary algebraic and geometric results}
\label{sec:prel-results-about}

We introduce several algebraic and algorithmic results about finite field
extensions and elliptic curves over finite fields

\subsection{Finite field extensions}
\label{sec:finite-field-extens}
In this section,
we collect algebraic and algorithmic results about finite fields.

\subsubsection{Generator of a subextension}\label{section:gene}
 
We prove the following lemma.

\begin{lemma}[Subfield generated by a symmetric function]\label{lemma:subfield}
  Let $\bM$ be a finite field and let $\bK$ be a subfield of $\bM$. We assume
  that the degree of $\bM$ over $\bK$ is a power of a prime integer $\ell$.
  Let $\alpha$ be a generator of $\bM$ over $\bK$.  Let $\bL$ be a subfield of
  $\bM$ containing $\bK$. Let $n$ be the degree of $\bM$ over $\bL$.  Let
  $(\Sigma_k(\alpha))_{0 < k \leqslant n}$ be the $n$ symmetric functions of
  $\alpha$ above $\bL$.  Then at least one among these $n$ symmetric functions
  generates $\bL$ over $\bK$.
\end{lemma}

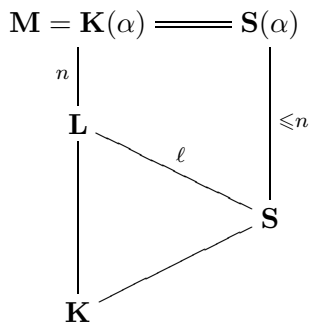
\begin{displaymath}
    \xymatrixrowsep{0.3in}
    \xymatrixcolsep{0.4in}
  \xymatrix{ \bM=\bK(\alpha) \ar@{=}[r] \ar@{-}[d]_{n} & \bS(\alpha)
    \ar@{-}[dd]^{\leqslant n} \\ \bL \ar@{-}[rd]^{\ell} \ar@{-}[dd] & \\ &\bS
    \ar@{-}[ld] \\ \bK& }
\end{displaymath}

\begin{proof} 
  If $\bL = \bK$, there is nothing to prove. When $\bL$ is a non
  trivial extension of $\bK$, the degree of this extension is a power of
  $\ell$.
  Let $\bS$ be the unique maximal proper subfield of $\bL$ containing $\bK$. The
  degree of $\bL$ over $\bS$ is $\ell$.
  The extension $\bM/\bL$ is cyclic of finite degree $n$, a power
of $\ell$.

  The $n$ elementary
symmetric functions of $\alpha$ over $\bL$ are the coefficients of
  the minimal polynomial of $\alpha$, seen as an element in the $\bL$-algebra
  $\bM$.
  We claim that at least one of these symmetric functions
  $(\Sigma_k(\alpha))_{0 < k \leqslant n}$ generates $\bL$ over $\bK$.
  Otherwise, all these functions would be contained in $\bS$.
  The field $\bS(\alpha)$ would then be a degree $\leqslant n$ algebraic
  extension of $\bS$.  Since $\bS(\alpha)$ contains $\bK(\alpha)$,
  $\bS(\alpha)$ is $\bM$.  But the degree of $\bM$ over $\bL$ is $n$, and this
  is greater than or equal to the degree of $\bM$ over $\bS$. So $\bL=\bS$. A
  contradiction.
\end{proof}

\subsubsection{Compositum}\label{section:general}

In this section, $\bK$ is a finite field with $q=p^w$ elements and $\Omega$ is
an algebraic closure of $\bK$. For every positive integer $k$, we denote by
$\FF_{p^k}$ the unique subfield of $\Omega$ with $p^k$ elements.  We have seen
in Section~\ref{section:qm1} how to construct an irreducible polynomial of
degree $d$ in $\bK[x]$ when $d$ is a prime power $\ell^\delta$ and $\ell$
divides $p(q-1)$. In Sections~\ref{section:fiber} and~\ref{section:basis}, we
shall treat the case when $d$ is a prime power $\ell^\delta$ and $\ell$ is
prime to $p(q-1)$.
The last problem to be considered is thus the following one: given two
irreducible polynomials $f_1(x)$ and $f_2(x)$ in $\bK[x]$ with coprime degrees
$d _1$ and $d_2$, construct a degree $d_1d_2$ irreducible polynomial.

\begin{lemma}[Composed sum of two polynomials]\label{lemma:compositum}
  There exists a deterministic algorithm that on input a finite field
  $\bK=(\ZZ/p\ZZ)[z]/(m(z))$ with $q$ elements, two irreducible
polynomials $f_1$ and
  $f_2$ in $\bK[x]$ of coprime degree $d_1$ and $d_2$, computes a degree $d_1d_2$
irreducible   polynomial in $\bK[x]$ at the expense of
  $(d_1d_2)^{1+\varepsilon(d_1d_2)}\times (\log q)^{1+\varepsilon(q)}$
elementary operations.
\end{lemma}
\begin{proof} 
  Let $\alpha_1\in \Omega$ be a root of $f_1(x)$.  Let $\alpha_2\in \Omega$ be
  a root of $f_2(x)$. We first show that $\alpha_1+ \alpha_2$ generates an
  extension of degree $d_1d_2$ of $\Fq$. Indeed, let $\Phi \in
  \Gal(\Omega/\Fq)$ be an automorphism that fixes $\alpha_1+\alpha_2$\,,
  \begin {displaymath}
    \Phi_{}(\alpha_1+ \alpha_2) = \alpha_1+ \alpha_2\,.
  \end {displaymath}
  One deduces that
  \begin{math}
    \Phi_{}(\alpha_1) - \alpha_1= \alpha_2- \Phi_{}(\alpha_2)
  \end{math} 
  is an element $\gamma$ of the intersection $\Fq$ of $\FF_{q^ {d_1}} $ and
  $\FF_{q^ {d_2}}$.  The order of $\Phi$ acting on $\FF_{q^ {d_1}} $ divides
  $d_1$. So $\Phi^{d_1}(\alpha_1)-\alpha_1=d_1\gamma=0$.  We prove in the same
  way that $d_2\gamma=0$. Since $d_1$ and $d_2$ are coprime we deduce that
  $\gamma=0$.  Thus $\Phi_{}$ acts trivially on $\FF_{q^ {d_1}} =
  \Fq(\alpha_1)$ and on $\FF_{q^ {d_2}} = \Fq(\alpha_2)$, therefore also on
  their compositum $\FF_{q^{d_1d_2}}$. So $\alpha_1+\alpha_2$ generates this
  compositum.
  Note that the same argument proves that $\alpha_1 \alpha_2$ generates
  $\FF_{q^{d_1d_2}}$.

  It is thus enough to compute the minimal polynomial of the sum, resp. the
  product, of $\alpha_1$ and $\alpha_2$. For this task, one may follow works
  by Bostan, Flajolet, Salvy and Schost~\cite{BoFlSaSc06}, based on algorithms
  for symmetric power sums due to Kaltofen and Pan~\cite{KP} and
  Sch\"onhage~\cite{Sch}.  The resulting polynomial is called the
  \textit{composed sum}, resp. the \textit{composed product}, of $f_1$ and
  $f_2$. See also \cite{G}. This yields an algorithm with a quasi-linear time
  complexity in $d_1d_2$.
\end{proof}

\subsubsection{Fast composition}\label{subsection:comp}

The following theorems were recently proven by Umans and Kedlaya~\cite{UK}.

\begin{theorem}[Kedlaya and Umans]\label{th:compo} 
  There exists a deterministic algorithm that on input a finite field
  $\bK=(\ZZ/p\ZZ)[z]/(m(z))$ with $q$ elements and three polynomials $f(x)$,
  $g(x)$ and $h(x)$ in $\bK[x]$ with degrees bounded by $d$, outputs the
  remainder $f(g(x)) \bmod h(x)$ at the expense of $d^{1+\varepsilon(d)}(\log
  q)^{1+\varepsilon(q)}$ elementary operations.
\end{theorem}

\begin{theorem}[Kedlaya and Umans]\label{th:mini} 
  There exists a deterministic algorithm that on input a finite field
  $\bK=(\ZZ/p\ZZ)[z]/(m(z))$ with $q$ elements, a degree $d$ irreducible monic
  polynomial $f(x)$ in $\bK[x]$, and a degree $\leqslant d-1$ polynomial
  $g(x)$ in $\bK[x]$,
 outputs the minimal polynomial
  $h(x)\in \bK[x]$ of the class $g(x)\bmod f(x)$ at the expense of $d^{1+\varepsilon(d)}(\log
  q)^{1+\varepsilon(q)}$ elementary operations.
\end{theorem}

The following corollary of Theorem~\ref{th:compo} is particularly useful.
\begin{corollary}\label{cor:compof}
  There exists a deterministic algorithm that on input a finite field $\bK$
  with $q$ elements given by a quotient $(\ZZ/p\ZZ)[z]/(m(z))$ and two
  rational fractions $F(x)$ and $G(x)$ in $\bK(x)$ with respective degrees
  $d_F$ and $d_G$, outputs the composition $F(G(x))=u(x)/v(x)$ where $u(x)$
  and $v(x)$ are coprime polynomials, at the expense of
  $(d_Fd_G)^{1+\varepsilon(d_Fd_G)}(\log q)^{1+\varepsilon(q)}$ elementary
  operations.
\end{corollary}

We first notice that the problem is trivial if one of the two fractions has
degree $1$. Composing $F$ and $G$ with rational linear fractions we may assume
that $F(0)=G(0)=0$. We compute the Taylor expansions at $0$ of either
fractions and we compose them using the algorithm in
Theorem~\ref{th:compo}. We recover the numerator $u(x)$ and denominator $v(x)$
of the corresponding fraction using the fast extended Euclid
algorithm~\cite[Chapter 11]{GG}.

\subsection{Elliptic curve over finite fields}
\label{sec:elliptic-curve-over}

We now state 
several known and useful facts about elliptic curves over finite
fields.

\subsubsection{Quotient isogenies}\label{section:noy}

Let $\bK$ be a finite field of characteristic $p$ and cardinality $q$. Let $E$
be an elliptic curve over $\bK$.  We denote by $\varphi_E : E \rightarrow E$
the degree $q$ Frobenius endomorphism of $E$. Let $t$ be the trace of
$\varphi_E$.  Let $\cO$ be the quotient ring $\ZZ[X]/(X^2-tX+q)$ and let
$\alpha$ be the class of $X$ in $\cO$.  Let $\epsilon_E : \cO\rightarrow
\End(E)$ be the ring homomorphism that maps $\alpha$ onto $\varphi_E$. We say
that $\epsilon_E$ is the {\it standard labeling} of $E$.

Let $S$ be a subset of $\cO$ containing an integer that is
prime to $p$. We define the
{\it kernel} of $S$ in $E$ to be the intersection of the kernels of all
endomorphisms $\epsilon_E(s)$ for $s\in S$. This is a finite {\'e}tale
subgroup of $E$. So, it is characterized by its set of geometric points.
We denote it $E[S]$. 

Now let $F$ be another elliptic curve over $\bK$ and let $\iota : E\rightarrow
F$ be an isogeny defined over $\bK$.  Let $\epsilon_F : \cO\rightarrow
\End(F)$ be the morphism of free $\ZZ$-modules that sends $1$ onto the
identity and $\alpha$ onto $\varphi_F$.  For any element $s$ in $\cO$, we have
\begin{equation}\label{eq:comm} 
  \iota\circ \epsilon_E(s) =\epsilon_F(s)\circ \iota.
\end{equation} 
Indeed, the identity above is true for $s=\alpha$ because $\iota$ is defined
over $\bK$.  It is evidently true also for $s=1$. Therefore it is true for all
$s$ in $\cO$ by linearity.
We deduce from Eq.~(\ref{eq:comm}) that $\epsilon_F$ is a ring homomorphism,
just as $\epsilon_E$.

Now let $G$ be a third elliptic curve over $\bK$. Let $\jmath : F\rightarrow
G$ be an isogeny defined over $\bK$.  We define $\epsilon_G : \cO\rightarrow
\End(G)$ as before.
Assume $\iota: E\rightarrow F$ is separable with kernel $E[S]$ where $S$ is a
subset of $\cO$ containing a prime to $p$ integer.  Assume $\jmath:
F\rightarrow G$ is separable with kernel $F[T]$ where $T$ is a subset of $\cO$
containing a prime to $p$ integer. Then the kernel of 
\begin{displaymath}
 \jmath\circ \iota :  \xymatrix{ E \ar@{->}[r]^{\iota} & F \ar@{->}[r]^{\jmath} & G}
\end{displaymath}
is $E[ST]$.
Indeed, both the kernel of $\jmath\circ \iota$ and $E[ST]$ are \'etale, so
they are characterized by their geometric points.
Now let $x$ be a point in the kernel of $\jmath\circ \iota$. Its image
$\iota(x)$ by $\iota$ lies in the kernel of $\jmath$. Therefore it is killed
by $T$: for any element $t$ of $T$ one has $\epsilon_F(t)(\iota(x))=0_F$.  So
$\iota(\epsilon_E(t)(x))=0_F$ and $\epsilon_E(t)(x)$ belongs in the kernel of
$\iota$. Thus it is killed by $S$: for any $s$ in $S$ we have
$\epsilon_E(s)(\epsilon_E(t)(x))=0_E$ or equivalently $\epsilon_E(st)(x)=0$.
Therefore $x$ lies in $E[ST]$.

Conversely, let $x$ be a point in $E[ST]$.  Let $t$ be an element in $T$. We
observe that $\epsilon_E(t)(x)$ is killed by $S$, so it belongs to the kernel
of $\iota$. Thus $\iota(\epsilon_E(t)(x))=\epsilon_F(t)(\iota(x))=0_F$. So
$\iota(x)$ is killed by $T$, therefore it belongs to the kernel of
$\jmath$. Thus $\jmath(\iota(x))=0_G$.
\medskip

\begin{lemma}[Composition of quotient isogenies]\label{lemma:ker} 
  Let $\bK$ be a finite field with characteristic $p$. Let $E$ be an elliptic
  curve over $\bK$. Let $t$ be the trace of the Frobenius endomorphism of $E$.
  Let $\cO$ be the quotient ring $\ZZ[X]/(X^2-tX+q)$ and let $\epsilon_E :
  \cO\rightarrow \End(E)$ be the standard labeling.  Let $S$ be a subset of
  $\cO$ containing a prime to $p$ integer and let $\iota : E\rightarrow F$ be
  the quotient by $E[S]$ isogeny.  Let $T$ be a subset of $\cO$ containing a
  prime to $p$ integer and let $\jmath : F\rightarrow G$ be the quotient by
  $F[T]$ isogeny.
  Then the kernel of $\jmath\circ \iota$ is $E[ST]$.
\end{lemma}

We note that a    general and  conceptual study of the action
of rings over group schemes was initiated by Serre \cite{serre},
Giraud \cite{giraud}, Waterhouse \cite{Water}, and Lenstra
\cite{lens}. In the case of elliptic curves one may use canonical lifts
and reduce to complex multiplication theory 
in characteristic zero. We prefer
a more self-contained and elementary approach.

\subsubsection{Density of elliptic curves with an $\ell$-torsion
point}\label{section:densite}

Let $\bK$ be a finite field with $q$ elements and let $\ell$ be a prime
integer.  Lenstra~\cite{Lenstra87} and Howe~\cite{Howe} give estimates for the
density of elliptic curves over $\bK$ whose number of $\bK$-rational points is
divisible by $\ell$.  In this section, we recall what these authors mean by
density and we explain why this density fits with the uniform density on
Weierstrass curves.

We call $\cE(\bK)$ the set of $\bK$-isomorphism classes of elliptic curves
over $\bK$.  The $\bK$-isomorphism class of a curve $E/\bK$ is denoted by
$[E]$.  One defines a measure on the finite set $\cE(\bK)$ in the following
way: the measure of a class $[E]$ is the inverse of the group of
$\bK$-automorphisms of $E$.  So the measure of a subset $S$ of $\cE(\bK)$ is
\begin{equation}\label{eq:mes}
  \mu_\cE(S)=\sum_{[E]\in S} \frac{1}{\#\Aut_\bK(E)}.
\end{equation}
Lenstra and Howe prove that the measure of the full set $\cE(\bK)$ is
$q$.\medskip

Now, let $\cW(\bK)$ be the set of Weierstrass elliptic curves over $\bK$,
\begin{equation}\label{eq:we}
  E/\bK~: Y^2Z+a_1XYZ+a_3YZ^2=X^3+a_2X^2Z+a_4XZ^2+a_6Z^3\,.
\end{equation}
We denote by $\mu_\cW$ the uniform measure on this set: the $\mu_\cW$-measure
of a subset of $\cW(\bK)$ is defined to be its cardinality.  This is a very
convenient measure. In order to pick a random Weierstrass curve according to
this measure, we just choose each coefficient $a_1$, $a_2$, $a_3$, $a_4$,
$a_6$ at random with the uniform probability in $\bK$ and we check that the
discriminant is non-zero (if it is zero we start again).

Let $\gamma : \cW(\bK)\rightarrow \cE(\bK)$ be the map that to every curve $E$
associates its isomorphism class $[E]$.  This is a surjection: every elliptic
curve over $\bK$ has a Weierstrass model over $\bK$.
Let $\bA(\bK)$ be the group of projective transforms of the form
\begin{displaymath} 
  (X:Y:Z)\mapsto (u^2X+rZ : u^3Y + su^2X + tZ:Z)
\end{displaymath} 
where $u\in \bK^*$ and $r, s, t \in \bK$. This group acts on the set
$\cW(\bK)$ of Weierstrass elliptic curves over $\bK$.  Two Weierstrass
elliptic curves over $\bK$ are isomorphic over $\bK$ if and only if they lie
in the same orbit for the action of $\bA(\bK)$. Further the group of
$\bK$-automorphisms of a Weierstrass elliptic curve is isomorphic to the
stabilizer of $E$ in $\bA(\bK)$.

So the orbit of a Weierstrass curve $E/\bK$ under the action of $\bA(\bK)$ is
the fiber $\gamma^{-1}([E])$ and the cardinality of this fiber is the quotient
\begin{math}
  {\#\bA(\bK)}\,/\,{\#\Aut_\bK(E)}.
\end{math}
Therefore, if $S$ is a subset of $\cE(\bK)$ and if $T$ is its preimage by
$\gamma$, then the measures of $S$ and $T$ are proportional,
\begin{displaymath} 
  \mu_\cW(T)= \# \bA(\bK)\times \mu_\cE(S)\text{ where }\#\bA(\bK)=(q-1)q^3.
\end{displaymath}
In particular, if we want to pick a random $\bK$-isomorphism class of elliptic
curve according to the measure $\mu_\cE$, it suffices to pick a random
Weierstrass elliptic curve according to the uniform measure $\mu_\cW$.

We now can state a special case of the main result in Howe's
paper~\cite{Howe}.
\begin{theorem}[Howe]\label{th:Howe} 
  Let $q$ be a prime power and let $\bK$ a field with $q$ elements. Let
  $\cE(\bK)$ be the set of $\bK$-isomorphism classes of elliptic curves over
  $\bK$. Let $\mu_\cE$ be the measure on this set defined by
  Eq.~(\ref{eq:mes}). Let $\ell$ be a prime integer not dividing $q-1$. The
  isomorphism classes in $\cE(\bK)$ of elliptic curves having a $\bK$-rational
  point of order $\ell$ form a subset of density $r(\ell,q)$ where
  \begin{displaymath}
    \left| r(\ell,q)-\frac{1}{\ell-1}\right| \le     \frac{4\ell(\ell+1)}{(\ell-1)\sqrt q}.
  \end{displaymath}
\end{theorem}

We deduce the following corollary.
\begin{corollary}[Density of elliptic curves with an $\ell$-torsion point]
  \label{cor:Howe}
  Let $q$ be a prime power and let $\bK$ be a field with $q$ elements. Let
  $\cW(\bK)$ be the set of Weierstrass elliptic curves over $\bK$. Let
  $\mu_\cW$ be the uniform measure on this set. Let $\ell$ be a prime integer
  not dividing $q-1$. The density  $r(\ell,q)$
of Weierstrass curves having a
  $\bK$-rational point of order $\ell$ satisfies 
  \begin{displaymath}
    \left| r(\ell,q)-\frac{1}{\ell-1}\right| \le     \frac{4\ell(\ell+1)}{(\ell-1)\sqrt q}.
  \end{displaymath}

\end{corollary}

\section{Isogeny fibers}\label{section:fiber}

In this section, we show how to construct irreducible polynomials using
elliptic curves.  

Let $\bK$ be a field and let $\Omega$ be an algebraic closure of $\bK$. Let
$E/\bK$ be an elliptic curve given by the Weierstrass equation~\eqref{eq:we}.
We denote by $O_E=[0:1:0]$ the origin of $E$ and by $x=X/Z$, $y=Y/Z$ the
affine coordinates associated with the projective coordinates $[X: Y: Z]$.

Let $E'/\bK$ be another elliptic curve in Weierstrass form. We define $X'$,
$Y'$, $Z'$, $a_1'$, $a_2'$, $a_3'$, $a_4'$, $a_6'$, $x'$, $y'$, $O'$
similarly.  Let $\iota /\bK: E/\bK\rightarrow E'/\bK$ be a degree $d$
separable isogeny. We assume that $d$ is a positive odd number and the kernel
$\Ker \iota$ is cyclic.  Let $T\in E(\Omega)$ be a generator of $\Ker\iota$.
Let $\psi_\iota(x)\in \bK[x]$ be the degree $(d-1)/2$ polynomial
\begin{equation}\label {eq:psiI} 
  \psi_\iota(x) = \prod_{ 1 \leqslant k \leqslant (d-1) /2}(x-x(kT))\,.
\end{equation}
There exists a degree $d$ polynomial $\phi_\iota(x)\in \bK[x]$ and a
polynomial $\omega_\iota (x,y)=\omega_0(x)+y\omega_1(x)$ in $\bK[x,y]$ with
degree $1$ in $y$ such that the image of the point $(x,y)$ by $\iota$ is
$(x',y')$ where $x'={\phi_\iota(x)}/{\psi_\iota^2(x)}$ and
$y'={\omega_\iota(x,y)}/{\psi_\iota^3(x)}$. We call $I(x)$ the
rational fraction $I(x)={\phi_\iota(x)}/{\psi_\iota^2(x)}$.

Now let $A$ be a $\bK$-rational point on $E'$ such that $2A\not =O'$ and let
$B \in E(\Omega)$ be a point on $E$ such that $\iota (B)=A$.  We define the
degree $d$ 
polynomial
\begin{equation*}
  f_{\iota,A}(x) = \phi_\iota(x)-x'(A)\psi_\iota^2(x)\in \bK[x].
\end{equation*}
Its roots are the $x(B+kT)$ for $0 \leqslant k < d$, and they are pairwise
distinct because $2A\not = O'$.  So $f_{\iota,A}(x)$ is a degree $d$ separable
polynomial.
The coordinate $x$ lies in the field $\bK(E)$ of $\bK$-rational
functions on $E$. So the map $x : E(\Omega)-O\rightarrow \Omega$ 
induces a Galois equivariant
bijection between the fiber $\iota^{-1}(A)$ and the roots of $f_{\iota,A}(x)$.
In particular, $f_{\iota,A}(x)$ is irreducible if and only if the fiber
$\iota^{-1}(A)$ is.
The fiber  $\iota^{-1}(A)$ over $A$ is an affine scheme with ring
 $\bK[x]/(f_{\iota,A}(x))$ and the class of
$y$ in this ring is given by 
\begin{equation}\label{eq:y}
  y=\frac{y'(A)\psi_\iota^3(x)-\omega_0(x)}{\omega_1(x)}\bmod f_{\iota,A}(x)\,.
\end{equation}
\medskip

Then, the two questions that worry us are the following ones.
\begin{enumerate}[itemsep=0pt,parsep=0pt,partopsep=0pt,topsep=0pt]
\item Can we compute $f_{\iota,A}(x)$ quickly, \textit{e.g.} in quasi-linear
time in $d$~?
\item Under which conditions is $f_{\iota,A}(x)$ irreducible~?
\end{enumerate}
\medskip

These two questions are successively addressed in
Sections~\ref{subsection:fast-calc-polyn} and~\ref{subsection:conditions}.  In
Section~\ref{subsection:existence} we deduce a fast algorithm that constructs
a degree $d$ irreducible polynomial in $\bK[x]$ when $\bK$ is a finite field
with $q=p^w$ elements and $d=\ell^\delta$ is a power of a prime $\ell$ such
that $\ell$ is prime to $p(q-1)$ and $4\ell \leqslant q^\frac{1}{4}$. Larger
primes $\ell$ are considered in Section~\ref{section:basis}.

\subsection{Calculation of the polynomial $f_{\iota,A}(x)$}
\label{subsection:fast-calc-polyn}

For any geometric point $P \in E(\Omega)$, we denote by $ \tau_P: E
\rightarrow E$ the translation by $P$. Let $x_P$ be the function $x \circ
\tau_{- P} $ and similarly let $y_P$ be the function $y \circ \tau_{- P}$.
If $P=kT$, we moreover define $x_k=x_{kT} $ and $y_k=y_{kT} $.  Recall $d$ is
assumed to be odd.

In this section, we present methods for fast construction of
isogenies. Section~\ref{subsubsection:velu} concerns isogenies with split
cyclic kernel and V{\'e}lu's formulae.
Section~\ref{subsubsection:composition} recalls how one can take advantage of
the decomposition of  an isogeny into several ones with smaller degrees. This
is particularly useful when $E/\bK$ has complex multiplication and the isogeny
in question is the quotient isogeny associated to some power of an invertible
prime ideal in the endomorphism ring of $E$. This idea is detailed in
Section~\ref{subsubsection:simple}.

\subsubsection{V{\'e}lu's isogenies}\label{subsubsection:velu}

Let $T$ be a $\bK$-rational point and let $\iota$ be the isogeny given by
V{\'e}lu's formulae~\cite{Velu71},
\begin {equation} \label {eq:xyp} 
  \left\{
    \begin{array}{rcl} 
      x'&=&\displaystyle x+ \sum_{0 < k < d} \left (x_{k} - x(kT) \right)\,,\\
      y'&=&\displaystyle y+ \sum_{0 < k < d} \left (y_{k} - y(kT) \right)\,.
    \end{array}\right.
\end {equation}
We put some order in Eq.~(\ref{eq:xyp}). Using the addition law on $E$, we
first express $x_k$ in terms of $x$ and $y$,
\begin {multline} \label {eq:xk} 
  x_{kT} \cdot (x-x(kT))^2
  = x(kT) x^2+ \left(a_3+2y(kT) +a_1x(kT) \right) y  \\ 
  + \left(a_4+a_1^2x(kT) +a_1a_3+2a_2x(kT) +a_1y(kT) +x(kT) ^2 \right) x  \\
  + a_3^2+a_1a_3x(kT) +a_3y(kT) +a_4x(kT) +2a_6.
\end{multline}
We deduce that $(x_{kT} +x_{- kT} - 2x(kT))\cdot (x-x(kT))^2$ is equal to
\begin {multline} \label {eq:xkxmk} (6x(kT) ^2+(a_1^2+4a_2) x(kT)
  +a_1a_3+2a_4)\, x -2x(kT) ^3\\
  +(a_1a_3+2a_4) x(kT) +a_3^2+4a_6\,.
\end {multline}
One computes the rational fraction $x'={\phi_\iota(x)}/{\psi_\iota^2(x)}$
using Eqs.~(\ref{eq:xyp}) and (\ref{eq:xkxmk}) by gathering the terms relative
to $k$ and $-k$, with the help of a ``divide and conquer''
strategy~\cite[Chapter 10]{GG}.  Complexity is quasi-linear in $d$.
A similar calculation gives us the explicit form of
$y'={\omega_\iota(x,y)}/{\psi_\iota^3(x)}$.
\medskip

\noindent {\it Example.}
We take $p=7$, $q=7$ and $d=5$. The elliptic curve
$E/\FF_7$ with equation 
${y}^{2}={x}^{3}+x+4$ has got ten $\FF_7$-rational points. The point
$T = (6, 4)$ has order $\ell=5$.  The group generated by $T$ is
\begin{equation*} 
  \langle T \rangle = \left\{ O_E,\ (6, 4),\ (4, 4),\ (4, 3),\ (6, 3) \right\}\,.
\end{equation*} 
The corresponding isogenous curve $E'$ is given by V\'elu's formulae,
$E':{y}^{2}={x}^{3}+3\,x+4.$
Moreover, Eq.~\eqref{eq:xyp} yields
\begin{multline*} 
  x' = x+ {\frac {y+6\,{x}^{2}+2\,x}{ \left( x+1 \right) ^{2}}} - 6 + {\frac
    {y+4\,{x}^{2}+3\,x+5}{ \left( x+3 \right) ^{2}}} - 4 +\\
  {\frac {6\,y+4\,{x}^{2}+3\,x+5}{ \left( x+3 \right) ^{2}}} - 4 + {\frac
    {6\,y+6\,{x}^{2}+2\,x}{ \left( x+1 \right) ^{2}}} - 6\,.
\end{multline*} 
Using Eq.~\eqref{eq:xkxmk}, we find an expression for $x'$ in terms of $x$
alone,
\begin{displaymath}
  x' = x + {\frac {x+2}{ \left( x+1 \right) ^{2}}} + \frac {1}{ \left( x+3
    \right) ^{2}} = 
  {\frac {{x}^{5}+{x}^{4}+2\,{x}^{3}+5\,{x}^{2}+4\,x+5}{ \left( x+3 \right) ^{2} \left(
        x+1 \right) ^{2}}}\,.
\end{displaymath} 
It remains to choose a point $A$  in $E'(\FF_7)$. We set 
$A=(1,1)$, a point of order $5$, and we finally obtain,
\begin{displaymath} 
  f_{\iota,A}(x) =
  {x}^{5}+{x}^{4}+2\,{x}^{3}+5\,{x}^{2}+4\,x+5 - { \left( x+3 \right) ^{2}
    \left( x+1 \right) ^{2}} = {x}^{5}+{x}^{3}+4\,{x}^{2}+x+3\,.
\end{displaymath}

\subsubsection{Composition of isogenies}\label{subsubsection:composition}

Assume $d$ factors as $d_1d_2$. Then the degree $d$ isogeny $\iota :
E\rightarrow E'$ decomposes as $\iota = \iota_2\circ \iota_1$ where $\iota_1 :
E\rightarrow F$ is a degree $d_1$ isogeny and $\iota_2 : F\rightarrow E_2$ is
a degree $d_2$ isogeny.  The kernel of $\iota_1$ is generated by $d_2T$ and
the kernel of $\iota_2$ is generated by $\iota_1(T)$.  Let $I(x)$ be the
degree $d$ rational fraction associated with $\iota$. Define similarly
$I_1(x)$ and $I_2(x)$. Then $I(x)=I_2(I_1(x))$.  We may then compute $I(x)$ in
three steps: first compute $I_1(x)$, then compute $I_2(x)$, and finally
compute the composition $I=I_2\circ I_1$ using work by Umans and
Kedlaya~\cite{Umans08,UK} (see Corollary~\ref{cor:compof} in
Section~\ref{subsection:comp}).

\subsubsection{A special simple case}\label{subsubsection:simple}

We now assume that $\bK$ is a finite field with $q=p^w$ elements.  Let
$\varphi_E : E \rightarrow E$ be the Frobenius endomorphism of $E$ and $t$ be
its trace. Let $\cO$ be the quotient ring $\ZZ[X]/(X^2-tX+q)$ and let $\alpha$
be the class of $X$ in $\cO$.  We call $\epsilon : \cO \rightarrow \End(E)$
the ring monomorphism that sends $\alpha$ onto $\varphi_E$. For every subset
$S$ of $\cO$, we define the {\it kernel} of $S$ in $E$ to be the intersection
of all the kernels of the endomorphisms $\epsilon(s)$ for $s\in S$. We denote
it by $E[S]$.
Let $\ell $ be a prime not dividing $p(q-1)$. We assume that $\ell$ divides
the order $q+1-t$ of $E(\bK)$.  As a consequence $\ell$ is coprime to the
discriminant $t^2-4q$ of $\cO$.

We have
\begin{displaymath}
  X^2-tX+q=(X-1)(X-q) \bmod \ell\,,
\end{displaymath}
because $1-t+q$ is divisible by $\ell$ and the product of the roots of
$X^2-tX+q$ equals $q$. Furthermore, the roots $(1 \bmod \ell)$ and $(q \bmod
\ell)$ are distinct because $\ell$ does not divide $q-1$.
Let $\lgot=(\ell, \alpha -1)$ be the prime ideal in $\cO$ above $\ell$ and
containing $\alpha -1$. The kernel of $\lgot$ in $E$ is $E[\ell](\bK)$, the
rational part of the $\ell$-torsion of $E$.  This is a cyclic group of order
$\ell$ because $\ell$ divides $q+1-t$ and $\ell$ is coprime to $p(q-1)$.

Let $k$ be a positive integer. According to Hensel's lemma, there exist two
integers $\lambda_k$ and $\mu_k$ in $[0,\ell^k[$ such that $\lambda_k=1\bmod
\ell$, $\mu_k=q\bmod \ell$ and
\begin{displaymath}
  X^2-tX+q=(X-\lambda_k) (X-\mu_k)\bmod \ell^k.
\end{displaymath}
The ideal $\lgot^k$ of $\cO$ is generated by $\ell^k$ and $\alpha-\lambda_k$.
We show that the kernel $E[\lgot^k]$ of $\lgot^k$ in $E$ (in the sense of
section~\ref{section:noy}) is a cyclic group of order $\ell^k$ inside
$E(\Omega)$.
Let $\iota_k : E\rightarrow E_k$ be the quotient isogeny by $E[\lgot^k]$.
The elliptic curve $E_k$ is defined over $\bK$, a finite field with $q$
elements.  Let $\epsilon_k : \cO\rightarrow \End(E_k)$ be the ring
homomorphism that sends $\alpha$ onto the $q$-Frobenius endomorphism of
$E_k$. The two homomorphisms $\epsilon$ and $\epsilon_k$ are compatible with
the isogeny $\iota_k$ in the sense that for every $s$ in $\cO$ one has
$\epsilon_k(s)=\iota_k\circ \epsilon(s) \circ \iota_k^{-1}$.
Using Lemma~\ref{lemma:ker} of Section~\ref{section:noy}, we see that
$\iota_{k+1} : E\rightarrow E_{k+1}$ decomposes as $\jmath_{k+1}\circ \iota_k$
where $\jmath_{k+1} : E_k\rightarrow E_{k+1}$ is the degree $\ell $ isogeny
with kernel $E_k[\lgot]=E_k[\ell](\bK)$. In particular, $\iota_k$ has degree
$\ell^k$, so the order of the kernel $E[\lgot^k]$ of $\iota_k$ is $\ell^k$.
This kernel is a subgroup of $E[\ell^k]$ that does not contain the full $\ell$
torsion, therefore it is cyclic.
We obtain in this way a chain of degree $\ell$ isogenies
\begin{displaymath}
  \xymatrixrowsep{0.3in}
  \xymatrixcolsep{0.3in}
\xymatrix{ E \ar@{->}[r]^{\jmath_1} & E_1\ar@{->}[r]  & \dots \ar@{->}[r] & E_k
  \ar@{->}[r]^{\jmath_{k+1}} & E_{k+1} \ar@{->}[r] & \dots}
\end{displaymath} 

We denote by $I_k(x)\in \bK(x)$ the degree $\ell^k$ rational fraction
associated with $\iota_k$. We denote by $J_k\in \bK(x)$ the degree $\ell$
rational fraction associated with $\jmath_k$.  In order to compute
$\jmath_1=\iota_1$, we pick a random point in $E(\bK)$ and multiply it by
$(q+1-t)/\ell$.  If the result is non-zero, we are done, otherwise we start
again.  We then compute  $I_1=J_1$ using V{\'e}lu's formulae in
Section~\ref{subsubsection:velu}. 
Every rational
fraction $J_k$ can be computed   the same way.
The composition $I_k=J_{k}\circ \cdots J_2\circ J_1$ can be computed using
the method in Paragraph~\ref{subsubsection:composition}.

\subsection {Irreducibility conditions} \label {subsection:conditions}

We assume that we still are in the situation of
Paragraph~\ref{subsubsection:simple}. Let $\ell $ be a prime not dividing
$p(q-1)$. In particular $\ell$ is odd. We assume that $\ell$ divides the order
$q+1-t$ of some elliptic curve $E/\bK$. 
As a consequence $\ell$ is coprime to $t^2-4q$.  We want to construct an
irreducible polynomial $f(x)\in \bK[x]$ with degree $d=\ell^\delta$. We factor
$q+1-t$ as $q+1-t=\ell^e\ell'$ where $\ell'$ is prime to $\ell$.

There exist two integers $\lambda_{e+\delta}$ and $\mu_{e+\delta}$ such that
\begin{eqnarray*} 
  \lambda_{e+\delta} = 1 \bmod \ell^{e}&,& \mu_{e+\delta} = q \bmod \ell^{e}\,, \\
  X^2-tX+q &=&(X - \lambda_{e+\delta})(X - \mu_{e+\delta}) \bmod \ell^{e+\delta}\,.
\end{eqnarray*}
We write $\lambda_{e+\delta} = 1 + \ell^e\ell''$ with $\ell''$ prime to
$\ell$. In the sequel, we set $\lambda=\lambda_{e+\delta}$ and
$\mu=\mu_{e+\delta}$.  Let now
\begin{displaymath}
  \dgot =(d, \alpha - \lambda) = (\ell, \alpha - \lambda) ^ {\delta} =\lgot^\delta.
\end{displaymath}
This is an invertible ideal. Its kernel $E[\dgot]$ in $E$ is the kernel of the
isogeny $\iota_\delta : E\rightarrow E_\delta$.  The $\ell$-Sylow subgroup of
$E_\delta(\bK)$ is the kernel of $\lgot^e=(\ell^e,\alpha-1)$ in $E_\delta$ and
it is cyclic. Let $A$ be a generator of it. Let $B \in E(\Omega)$ such that
$\iota_\delta(B)=A$.
Then, $B$ generates the kernel of $\lgot^{e+\delta}= (\ell^{e+\delta},
\varphi_E - \lambda)$ in $E$.
Especially,
\begin{equation}\label{eq:lambda}
  \varphi_E (B) = \lambda B,
\end{equation} 
and the order of $\lambda=1+\ell^e\ell''$ in $(\ZZ/\ell^{e+\delta}\ZZ)^*$ is
$d=\ell^\delta$.  Thus, the Galois orbit of $B$ has cardinality $d$ and the
polynomial $f_{\iota,A}(X)$ is irreducible.

\subsection{Existence conditions}\label{subsection:existence}

Assume we are given a finite field $\bK$ with characteristic $p$ and
cardinality $q$ and an integer $d=\ell^\delta$ such that $\ell$ is prime to
$p(q-1)$. We look for a degree $d$ irreducible polynomial in $\bK[x]$.  The
construction in Section~\ref{subsection:conditions} requires an elliptic curve
over $\bK$ such that $\ell$ divides the cardinality $q+1-t$ of $E(\bK)$.  Is
there any such elliptic curve ? How can we find it ?
\medskip

If $\ell\leqslant 2\sqrt q$, then there are at least two consecutive integer
multiples of $\ell$ in the interval $[q+1-2\sqrt q,q+1+2\sqrt q]$. At least
one of them is not congruent to $1$ modulo $p$.  So there exists at least one
elliptic curve with cardinality divisible by the prime $\ell$.
We want to bound from below the number of such elliptic curves.  We use the
results of Lenstra~\cite{Lenstra87} extended by Howe~\cite{Howe}.

From Theorem~\ref{th:Howe} and Corollary~\ref{cor:Howe} of
Section~\ref{section:densite}, we deduce that the probability that a
Weierstrass elliptic curves over a finite field $\bK$ with $q$ element has an
order divisible by $\ell$ is ${1}/(\ell -1)$, up to an error term bounded in
absolute value by ${8\ell}/{\sqrt q}$. We deduce that if
\begin{equation}\label{eq:cond}
  4\ell\leqslant {q^{\frac{1}{4}}}{},
\end{equation}
then this proportion is at least ${1}/({2\ell})$.

In that case, we can find such an elliptic curve in the following way: we pick
a random Weierstrass elliptic curve over $\bK$. We compute its cardinality
using Schoof's algorithm at the expense of $(\log q)^{5+\varepsilon(q)}$
elementary operations.  If this cardinality is divisible by $\ell$ we are
done. Otherwise, we try again. The average number of trials is $O(\ell )$.
The expected time to find the needed curve $E$ is $\ell (\log
q)^{5+\varepsilon(q)}$ elementary operations provided condition
(\ref{eq:cond}) holds true.
\medskip

All in all, we need $\ell \times (\log q)^{5+\varepsilon(q)}$
elementary operations to find the first elliptic curve, then
$\delta^{1+\epsilon(\delta)} \times
  \ell^{1+\varepsilon(\ell)} \times (\log q)^{2+\varepsilon(q)}$
elementary operations to compute the $\delta$ isogenies
of degree  $\ell$, and 
$d^{1+\varepsilon(d)}\times (\log q)^{1+\varepsilon(q)}$ elementary
operations to compose these isogenies.
The conclusion of this section is the following.

\begin{lemma}[Isogeny fiber]\label{lemma:fiber}
  There exists a probabilistic (Las Vegas) algorithm that on input a finite
  field $\bK$ with characteristic $p$ and cardinality $q=p^w$, a prime integer
  $\ell$ not dividing $p(q-1)$ such that $4\ell\leqslant {q^{\frac{1}{4}}}{}$,
  and a positive integer $\delta$, computes an irreducible polynomial in
  $\bK[x]$ of degree $d = \ell^\delta$, at the expense of 
$\ell \times (\log q)^{5+\varepsilon(q)} + d^{1+\varepsilon(d)}\times (\log q)^{2+\varepsilon(q)}$ elementary
  operations.
\end{lemma}

\section{Base change}\label{section:basis}

In this section, $\bK= (\ZZ/p\ZZ)[z]/(m(z))$ is a finite field with $q=p^w$
elements.  We still assume here that $d=\ell^\delta$ is a power of a prime
$\ell$ where $\ell$ is prime to $p(q-1)$. We look for a degree $d$ irreducible
polynomial in $\bK[x]$.  However, we no longer assume that $4\ell \leqslant
q^{\frac{1}{4}}$.

We adapt the main idea in Section~\ref{subsection:gauss} to the context of
elliptic curves: we base change to a small auxiliary extension.
Let $n$ be the smallest integer coprime with $\ell (\ell-1)$ such that $Q=q^n$
satisfies $4\ell \leqslant Q^{\frac{1}{4}}$.  According to Iwaniec's result
about Jacobsthal's problem~\cite{Iwa} we have $n=(\log \ell) ^{2 +
\varepsilon(\ell)}$.  Let us remark that $d$ is then coprime with $Q-1$ too.

Using e.g. the methods in Shoup~\cite{Shoup93}, we find a degree $n$
irreducible polynomial $m'\in \bK[z']$. We set $\bL=\bK[z']/(m'(z'))$. A basis
of this $(\ZZ/p\ZZ)$-vector space is given by the $z^jz'^i$ for $0\leqslant i
<n$ and $0\leqslant j < w$.  Using the method given in the introduction, we
find a generator $\tau$ of the $(\ZZ/p\ZZ)$-algebra $\bL$.  We compute also
the minimal polynomial $h(u)\in (\ZZ/p\ZZ)[u]$ of $\tau$. We set $\tilde
\bL=(\ZZ/p\ZZ)[u]/(h(u))$.  A basis of this $(\ZZ/p\ZZ)$-vector space is given
by the $u^k$ for $0\leqslant k < nw$. We compute and store the matrix of the
isomorphism $\kappa : \tilde \bL \rightarrow \bL$ that sends $u\bmod h(u)$
onto $\tau$. This is a $nw \times nw$ matrix with entries in $\ZZ/p\ZZ$.  We
also compute and store the inverse of this matrix. The image $\tilde
\bK=\kappa^{-1}(\bK)$ of $\bK$ by $\kappa^{-1}$ is the unique subfield with
$q$ elements inside $\tilde\bL$.

The reason for introducing these two different models of the field with $q^n$
elements is that, on the one hand, this field should be constructed as an
extension of $\bK$ because we shall have to descend to $\bK$ later on; but on
the other hand, the field with $q^n$ elements should be also presented as a
monogenic extension of $\ZZ/p\ZZ$, because all the algorithms described and
used so far (and in particular the algorithms due to Umans and Kedlaya)
require that the base field be presented as a monogenic extension of
$\ZZ/p\ZZ$.
One can now apply the construction of Section~\ref{section:fiber} to $\tilde
\bL$ and obtain an irreducible polynomial $F_{\iota,A}(x)$ of degree $d$ in
$\tilde \bL[x]$, in time
\begin{math}
  (\log Q)^ {5+\varepsilon(Q)}d^{1+\varepsilon(d)},
\end{math}
that is
\begin{displaymath} 
  (\log q)^ {5+\varepsilon(q)}d^{1+\varepsilon(d)}\,
\end{displaymath}
elementary operations.
 
It remains to derive from $F_{\iota,A}(x)$ an irreducible polynomial $f(x)$ of
degree $d$ over $\bK$.
Recall $F_{\iota,A}(x)$ is the minimal polynomial of $x(B)$ where $B$ is a
geometric point of order $\ell^{e+\delta}\leqslant 4Qd$ on an elliptic curve
$E$ over $\tilde \bL$.  We also are given an integer $\lambda$ such that
$0\leqslant \lambda <\ell^{e+\delta}$ and
\begin{equation}\label{eq:lambdaQ}
  \varphi_E (B) = \lambda B
\end{equation}
where $\varphi_E$ is the degree $Q$ Frobenius endomorphism of $E/\tilde\bL$.
We set $\alpha=x(B) \in \tilde \bL[x]/(F_{\iota,A}(x))$.  Recall $\Phi_q$ is
the map which raises to the $q$-th power. We have $\Phi_Q= \Phi_q^n$.
The field $\tilde\bL[x]/(F_{\iota,A}(x)) = \tilde\bL(\alpha)$ is a
degree $d$  extension
 of $\tilde\bL$.  For any integer $k$ between $1$ and $n$, one
denotes by $\Sigma_k(\alpha)$ the $k$-th symmetric function of the conjugates
of $\alpha$ over the subfield with $q^d$ elements. These conjugates
are 
\begin{math} 
  \alpha,\ \Phi_q^d(\alpha),\ \Phi_q^ {2d}(\alpha),\ \ldots,\ \Phi_q^{(n-1)
    d}(\alpha)\,.
\end{math}
Since $d$ is a {\it prime power}, we deduce from Lemma~\ref{lemma:subfield} of
Section~\ref{section:gene} that at least one among these $n$ symmetric
functions generates the degree $d$ extension of $\tilde\bK$. In other words,
there exists a $k$ between $1$ and $n$ such that the polynomial
\begin{displaymath} 
  \tilde f(x) = \prod_{0 \leqslant l < d}(x - \Phi_q^l(\Sigma_k(\alpha)))
\end{displaymath} 
is irreducible of degree $d$ in $\tilde\bK[x]\subset \tilde\bL[x]$\,.
\medskip

Three questions now worry us, that we consider in turn in
Sections~\ref{subsection:Sigma},~\ref{subsection:how-find-integer}
and~\ref{subsection:how-calc-char}.
\begin{enumerate}[itemsep=0pt,parsep=0pt,partopsep=0pt,topsep=0pt]
\item How to compute $\Sigma_k(\alpha)$ and its conjugates~?
\item How to find the good integer $k$~?
\item How to compute $\tilde f(x)\in \tilde\bK[x]$ starting from
$F_{\iota,A}(x)\in \tilde\bL[x]$~?
\end{enumerate}

\subsection {Computing symmetric functions}
\label{subsection:Sigma}

First, we compute $\beta=y(B) \in \tilde\bL[x]/(F_{\iota,A}(x))$ using
Eq.~(\ref{eq:y}).
Let now $l$ be an integer between $0$ and $dn-1$. We explain how to compute
$\alpha_l = \Phi_q^l(\alpha)$.  We write $l=r+ns$ with $0 \leqslant r < n$ and
$0 \leqslant s < d$. Then,
\begin{displaymath}
  \alpha_l=\Phi_q^l(\alpha) = \Phi_q^r(\Phi_Q^s(\alpha))\,.
\end{displaymath}

We first compute $\Phi_Q^s(\alpha) =x(\varphi_E^{s}(B))=x(\lambda^s B)$ using
Eq.~(\ref{eq:lambdaQ}). To this end, we write $\lambda^s =R\bmod
\ell^{e+\delta}$ where $0\leqslant R < \ell^{e+\delta}$ and we multiply the
$\ell^{\delta+e}$-torsion point $B\in E( \tilde\bL[x]/(F_{\iota,A}(x)) )$ by
$R$ using fast exponentiation.  This is done at the expense of a constant
times $(\log Q +\log d)$ operations in $ \tilde\bL[x]/(F_{\iota,A}(x)) $.
One then raises $\Phi_Q^s(\alpha)$ to the $q^r$-th power at the expense of at
most $n \log q$ operations modulo $F_{\iota ,A}(x)$. Thus, each conjugate is
computed at the expense of
\begin{math}
  d^ {1+\varepsilon(d)}(\log q) ^ {2+\varepsilon(q)}
\end{math} elementary operations.

To compute all the $(\Sigma_k(\alpha))_{0 < k \leqslant n}$, one computes the
$n$ conjugates $\alpha$, $\Phi_q^d(\alpha)$, \ldots, $\Phi_q^ {(n-1)
d}(\alpha)$ and one forms the corresponding polynomial of degree $n$.
Altogether, the computation of the symmetric functions
$(\Sigma_k(\alpha))_{0<k\leqslant n}$ requires
\begin{displaymath}
  d^ {1+\varepsilon(d)}(\log q)^{2+\varepsilon(q)}
\end{displaymath}
elementary operations.

\subsection{Finding a generating symmetric function}
\label{subsection:how-find-integer}

One seeks an integer $k$ between $1$ and $n$ such that $\Sigma_k(\alpha)$
generates an extension of degree $d$ of $\tilde\bK$ (there is at least one
such integer). So we successively test all the $k$ between $1$ and $n$.  As
$n$ is small, this is not a problem.
We know that $\Sigma_k(\alpha)$ generates the degree $d$ extension of
$\tilde\bK$ if and only if
\begin{displaymath} 
  \Phi_q^{\ell^{\delta-1}}(\Sigma_k(\alpha))\not = \Sigma_k(\alpha)\,,
\end{displaymath}
where $\ell^{\delta-1}$ is the unique maximal divisor of $d$.  This condition
is equivalent to
\begin{displaymath} 
  \Sigma_k(\Phi_q^{\ell^{\delta-1}}(\alpha))\not = \Sigma_k (\alpha),
\text{ or }
  \Sigma_k(\alpha_{\ell^{\delta-1}}) \not = \Sigma_k(\alpha).
\end{displaymath}

One computes the $\Sigma_k(\alpha_{\ell^{\delta-1}})$'s in the same way as the
$\Sigma_k(\alpha)$'s, following Section~\ref{subsection:Sigma}.
It is then easy to compare $\Sigma_k(\alpha_{\ell^{\delta-1}})$ and
$\Sigma_k(\alpha)$.
One can thus find $k$ in
\begin{displaymath}
  d^ {1+\varepsilon(d)}(\log q) ^{2+\varepsilon(q)}
\end{displaymath}
elementary operations.

\subsection{Computing minimal polynomials}
\label{subsection:how-calc-char}

We now have an element $\Sigma_k(\alpha)$ of $\tilde\bL[x]/(F_{\iota,A}(x))$
and we know that it actually belongs to the degree $d$ extension of
$\tilde\bK$. But this is not really visible because $\Sigma_k(\alpha)$ is
given in the basis $1$, $x$, \ldots, $x^ {d-1} $ of $\tilde\bL [x]
/(F_{\iota,A}(x))$.
Still, the minimal polynomial $\tilde f(x)$ of $\Sigma_k(\alpha)$ has
coefficients in $\tilde\bK \subset \tilde\bL$. We compute this minimal
polynomial.  We use a general algorithm for this task, such as the one
appearing in recent work by Umans and Kedlaya~\cite{Umans08,UK}.  See
Theorem~\ref{th:mini} in Section~\ref{subsection:comp}. This algorithm
requires $d^{1+\varepsilon(d)}\times (\log Q)^{1+\varepsilon(Q)}$ elementary
operations. Finally, we apply the isomorphism $\kappa : \tilde\bL\rightarrow
\bL$ to every coefficient in $\tilde f(x)$ and we find a polynomial $f(x)$
with coefficients in $\bK\subset \bL$. This polynomial is irreducible in
$\bK[x]$.
\medskip

All in all, the conclusion of this section is the following.
\begin{lemma}[Base change]\label{lemma:main} 
  There exists a probabilistic (Las Vegas) algorithm that on input a finite
  field $\bK$ with characteristic $p$ and cardinality $q=p^w$, a prime integer
  $\ell$ not dividing $p(q-1)$, and a positive integer $\delta$, computes an
  irreducible polynomial in $\bK[x]$ of degree $d = \ell^\delta$, at the
  expense of $d^{1+\varepsilon(d)}\times (\log q)^{5+\varepsilon(q)}$
  elementary operations.
\end{lemma}

\section{Construction of irreducible polynomials}\label{section:compter}

\subsection{Finding one irreducible polynomial}\label{section:algo}

Given that we represent the finite field $\bK$ in a reasonable way, as
explained in the introduction, we can now state our main result.
\begin{theorem}\label{th:main} 
  There exists a probabilistic (Las Vegas) algorithm that on input a finite
  field $\bK$ with characteristic $p$ and cardinality $q=p^w$, and a positive
  integer $d$, returns a degree $d$ irreducible polynomial in $\bK[x]$. The
  algorithm requires $d^{1+\varepsilon(d)}\times (\log q)^{5+\varepsilon(q)}$
  elementary operations.
\end{theorem}
\begin{proof}
  The algorithms runs as follows.

  We first factor the degree $d$ as $d=\prod_{i}\ell_i^{\delta_i}$, this
  requires $O(d)$ elementary operations. Then, Lemma~\ref{lemma:compositum}
  shows that it suffices to find an irreducible polynomial of degree
  $\ell_i^{\delta_i}$ for every $i$.  So we may assume that $d=\ell^\delta$ is
  a prime power.

  Then, we use the construction of Lemma~\ref{lemma:ASE} if $\ell=p$, and
  Lemma~\ref{lemma:KE} if $\ell$ divides $q-1$. Otherwise
  Lemma~\ref{lemma:fiber} applies.
\end{proof}

\noindent \textit{Remark.} One may use a faster algorithm to compute the
cardinality of elliptic curves, for instance the SEA algorithm, and hope to
gain a $\log q$ speedup. But, at the time of writing, it is not clear how the
probability of failure of the SEA algorithm can be rigorously related to the
Las Vegas behavior of our construction, and we finally prefer to state a
complexity result based on Schoof's algorithm only.

\subsection{Constructing  random irreducible polynomials}

Let $\bK$ be a finite field with cardinality $q$ and characteristic $p$.  Let
$d\geqslant 2$ be an integer. 
We just explained  how to quickly compute a degree $d$
irreducible polynomial in $\bK[x]$.
We stress that the polynomials generated by our algorithm
have a very special form. This might be a problem for 
some applications. In this section we explain how to construct
{\em random} polynomials.
We  need an estimate for  the number of degree $d$
irreducible monic polynomials in $\bK[x]$.  We recall and prove a very
classical lower bound~\cite[Ex. 3.26 and 3.27, page 142]{Lidl}.

\begin{lemma}\label{lemma:lidl}
  Let $\bK$ be a finite field with $q$ elements. Let $d\geqslant 2$ be an
  integer.  The density of irreducible polynomials among the degree $d$ monic
  polynomials is greater than or equal to
  \begin{displaymath} 
    \frac{1}{d}(1-\frac{q}{q-1}(q^{-\frac{d}{2}}-q^{-d})).
  \end{displaymath} 
  Let $\bL$ be a degree $d$ field extension of $\bK$.  The density of
  generators of the $\bK$-algebra $\bL$ is greater than or equal to
  \begin{displaymath}
    1-\frac{q}{q-1}(q^{-\frac{d}{2}}-q^{-d}).
  \end{displaymath}
\end{lemma} 

\begin{proof}
Let $\Omega$ be an algebraic closure of $\bK$ and let $\bL$ be the unique
degree $d$ extension of $\bK$ inside $\Omega$. Call $\cG_d$ the set of
generators of the $\bK$-algebra $\bL$. This is the set of all $\alpha$ in
$\bL$ such that $\bK(\alpha)=\bL$.  Let $\cI_d$ be the set of degree $d$
monic irreducible polynomials in $\bK[x]$. Let $\rho : \cG_d\rightarrow
\cI_d$ be the map that to every generator $\alpha$ associates its minimal
polynomial. Every polynomial $P(x)$ in $\cI_d$ has exactly $d$ preimages by
$\rho$, namely its $d$ roots.

To enumerate the degree $d$ monic irreducible polynomials, we just count the
generators of $\bL$ over $\bK$.  Let $\alpha$ be an element in $\bL$. If
$\alpha$ does not generate $\bL$, then it belongs to a smaller extension of
$\bK$ inside $\bL$.  Therefore the complementary set of $\cG_d$ in $\bL$ is
the union of all proper subfields of $\bL$ containing $\bK$.  These subfields
are in correspondence with the strict divisors of $d$.  To any such divisor
$D$, we associate the unique extension of $\bK$ with degree $D$.  It has $q^D$
elements. The set of strict divisors of $d$ is a subset of $\{1,2,3,4, \ldots,
\lfloor \frac{d}{2}\rfloor\}$.  So the number of elements in $\bL$ that do not
generate it over $\bK$ is upper bounded by
\begin{displaymath} 
  q+q^2+q^3+q^4+\cdots+q^{\lfloor \frac{d}{2}\rfloor} =
  q\frac{q^{\lfloor \frac{d}{2}\rfloor}-1}{q-1}\leqslant
  \frac{q}{q-1}(q^{d/2}-1).
\end{displaymath} 
The cardinality of $\cG_d$ is thus $\geqslant q^d-\frac{q}{q-1}(q^{d/2}-1)$
and the cardinality of $\cI_d$ is
\begin{displaymath}
  \geqslant \frac{q^d}{d}-\frac{q}{d(q-1)}(q^{d/2}-1). \qedhere
\end{displaymath}
\end{proof}

If $d\geqslant 2$, we deduce from Lemma~\ref{lemma:lidl}
that the  density of generators
is $\geqslant
1-{1}/({q-1})=({q-2})/({q-1})$. So $\geqslant {1}/{2}$ if $q\geqslant
3$.
If $q=2$ and $d\geqslant 4$ then this density is $\geqslant 1-2\times
2^{-2}={1}/{2}$.
If $q=2$ and $d$ equals $2$ (resp. $3$) then this density is ${1}/{2}$
(resp. ${3}/{4}$).
If $d=1$ then this density is $1$.
We deduce the following lemma.
\begin{lemma}[Density of generators]\label{lemma:lidl2} 
  Let $\bK$ be a finite field with $q$ elements. Let $d\geqslant 1$ be an
  integer.  The density of irreducible polynomials among the degree $d$ monic
  polynomials is greater than or equal to ${1}/{2d}$.

Let $\bL$ be a degree $d$ field extension of $\bK$.  The density of generators
in the $\bK$-algebra $\bL$ is greater than or equal to ${1}/{2}$.
\end{lemma}

As a corollary of Lemma~\ref{lemma:lidl2}, given an irreducible polynomial
$f(x)$ of degree $d$ computed with Theorem~\ref{th:main}, one can compute a
new completely \textit{random} irreducible polynomial $g(x)$ 
at the expense of
only $d^{1+\varepsilon(d)}\times (\log q)^{1+\varepsilon(q)}$ elementary
operations. Indeed, we choose a random element in $\bL$, the degree $d$ extension of
$\bK$ constructed from $f(x)$, and we use Theorem~\ref{th:mini} to compute
its minimal polynomial. 
We obtain an irreducible
polynomial $g(x)$ that has degree $d$
with probability greater than
$1/2$. So, we have a Las Vegas quasi-linear algorithm.

\newpage

\end{document}